\documentclass [a4paper]{article}
\usepackage{amsmath}
\usepackage{mathrsfs}
\usepackage{amsfonts}
\usepackage {graphicx}
\usepackage{color,xcolor}
\usepackage{amsmath,amssymb}
\input{epsf.sty}
\usepackage{epsfig}

\usepackage{amsthm}
\newtheorem{Cor}{Corollary}[section]
\newtheorem{rem}{Remark}[section]
\newtheorem{theorem}{Theorem}[section]
\newtheorem{lemma}{Lemma}[section]
\newtheorem{algorithm}{Algorithm}
\newtheorem{definition}{Definition}[section]

\def\be{\begin{equation}}
\def\ee{\end{equation}}
\def\br{\begin{eqnarray*}}
\def\er{\end{eqnarray*}}

\begin{document}
\begin{center}
\Large\bf{Inexact Shift-and-Invert Arnoldi for Toeplitz Matrix Exponential
}\\
%
\quad\\

\normalsize~Ting-ting Feng\footnote[1]{School of Mathematics and statistics, Jiangsu Normal University, Xuzhou, 221116, Jiangsu, P.R. China.
Email: {\tt tofengtingting@163.com}. This author is supported by the Postgraduate Innovation Project of Jiangsu Province under grant CXLX13\_968.},
~
Gang Wu\footnote[2]{Corresponding author (G. Wu). Department of Mathematics,
China University of Mining and Technology \& School of Mathematics and statistics, Jiangsu Normal University, Xuzhou, 221116, Jiangsu, P.R. China.
E-mail: {\tt gangwu76@126.com} and {\tt wugangzy@gmail.com}. This author is
supported by the National Science Foundation of China under grant 11371176, the Natural Science Foundation of Jiangsu Province under grant BK20131126, the 333 Project of Jiangsu Province, as well as the Talent Introduction Program of China University of Mining and Technology.},
~
Yimin Wei\footnote[3]{School of Mathematical Sciences and Shanghai Key Laboratory of
Contemporary Applied Mathematics, Fudan University, Shanghai, 200433,
P.R. China.
Email: {\tt ymwei@fudan.edu.cn}. This author is supported by the National Natural Science Foundation of China under
grant 11271084.}

\end{center}

\begin{abstract}
We revisit the shift-and-invert Arnoldi method
proposed in [S. Lee, H. Pang, and H. Sun. {\it Shift-invert Arnoldi
approximation to the Toeplitz matrix exponential}, SIAM J. Sci.
Comput., 32: 774--792, 2010] for numerical approximation to the
product of Toeplitz matrix exponential with a vector. In this
approach, one has to solve two large scale Toeplitz linear systems
in advance. However, if the desired accuracy is high, the cost
will be prohibitive. Therefore, it is interesting to investigate how
to solve the Toeplitz systems inexactly in this method. The
contribution of this paper is in three regards. First, we give a new
stability analysis on the Gohberg-Semencul formula (GSF) and define the GSF condition number of a Toeplitz matrix. It is
shown that, when the size of the Toeplitz matrix is large, our result is sharper than the one given in [M.
Gutknecht and M. Hochbruck. {\it The stability of inversion formulas
for Toeplitz matrices}, Linear Algebra Appl., 223/224: 307--324,
1995]. Second, we
establish a relation between the error of Toeplitz systems and the
residual of Toeplitz matrix exponential. We show that if the GSF
condition number of the Toeplitz matrix is medium sized, then the
Toeplitz systems can be solved in a low accuracy. Third, based on this
relationship, we present a practical stopping criterion for relaxing
the accuracy of the Toeplitz systems, and propose an inexact
shift-and-invert Arnoldi algorithm for the Toeplitz matrix
exponential problem. Numerical experiments illustrate the numerical
behavior of the new algorithm, and show the effectiveness of our
theoretical results.\\
{\mbox {\bf Keywords:}} Toeplitz matrix, Matrix exponential, Shift-and-invert Arnoldi, Gohberg-Semencul formula (GSF), GSF condition number.
\end{abstract}


\section{Introduction}

\setcounter{equation}{0}


Toeplitz matrices occur in a variety of applications in mathematics
and engineering such as complex and harmonic analysis, statistics,
signal and image processing, information theory, numerical analysis,
see \cite{CJ,CN,Ng2} and the references therein. In this
paper, we are interested in numerical approximation to the product
of Toeplitz matrix exponential with a vector
\begin{equation}\label{eqn11}
{\bf y}(t)=\exp(-tA){\bf v},
\end{equation}
where $t$ is a scalar, ${\bf v}$ is a given vector, and $-tA$ is a
real $n\times n$ large Toeplitz matrix whose spectrum is located in the left half plane. This
problem plays an important role in various application fields such
as computational finance \cite{LP,DAM}, numerical solution of
Volterra-Wiener-Hopf equations \cite{MA}, calculating the
Wiener-Hopf integral equations \cite{GHK}, and so on.

The Krylov subspace method is an efficient approach to approximate
the matrix exponential with a vector, especially when the matrix is very
large \cite{Higham1,ML,Saad}. Indeed, it is the twentieth
dubious way to compute the matrix exponential \cite{ML}. In this
type of method, the matrix is first projected into a much smaller
subspace, then the exponential is applied to the projected matrix, and finally the approximation is projected
back to the original large space \cite{Higham1,ML,Saad}. This procedure can be achieved
by using the Lanczos process for symmetric matrices or by the
Arnoldi process for non-symmetric matrices, while both procedures
require only matrix-vector multiplications.

The shift-and-invert
Arnoldi and Lanczos methods were widely investigated to speed up the
Arnoldi and the Lanczos methods for matrix exponential \cite{MN,JH}.
Recently, by making use of the shift-and-invert Arnoldi technique,
Toeplitz structure and the famous Gohberg-Semencul formula (GSF)
\cite{GSF}, Pang {\it et al.} proposed a shift-and-invert Arnoldi
method for Toeplitz matrix exponential \cite{LP,Pang}. An advantage of this approach is that it is
unnecessary to explicitly form or store the Toeplitz matrix and its
inverse, and each matrix-vector product can be realized in several Fast Fourier Transformations
(FFTs) \cite{LP,Pang}. In the first step of this approach, one has to
solve two large scale (non-Hermitian) Toeplitz linear systems in a
desired accuracy. However, if the desired accuracy is very high, the
cost for solving the Toeplitz linear systems will be very large,
especially for some ill-conditioned problems. Thus, it is
interesting to investigate how to solve the Toeplitz linear systems
inexactly in the shift-and-invert Arnoldi method for matrix
exponential.

In this paper, we first give a new stability analysis on the
Gohberg-Semencul formula in terms of 1-norm and 2-norm, and define the ``GSF condition number" of a Toeplitz matrix. It is shown that our
results are sharper than the one given in \cite{MHG} when the
Toeplitz matrix is large.
We then establish a relation between the
error of Toeplitz systems and the residual of Toeplitz matrix
exponential. Based on the relationship, we present a practical
stopping criterion for solving the Toeplitz systems inexactly.

This paper is organized as follows. In Section 2, we briefly introduce the shift-and-invert Arnoldi method for Toeplitz matrix exponential \cite{LP}.
In Section 3, we give a stability analysis on the Gohberg-Semencul formula and propose an inexact shift-and-invert Arnoldi algorithm.
Numerical results given in Section 4 show the efficiency of our new algorithm and the effectiveness of the theoretical results.

\section{The shift-and-invert Arnoldi method for Toeplitz matrix exponential}

\setcounter{equation}{0}
In the shift-and-invert
Arnoldi/Lanczos method \cite{RRA,LP,MN,Pang,JH}, the Krylov subspace is constructed
by using the matrix $(I+\gamma A)^{-1}$, where $\gamma$ is a user-prescribed parameter and $I$ is the identity matrix whose order is clear from context.
Let ${\bf v}_1={\bf v}/\|{\bf v}\|_2$, the $m$-step shift-and-invert Arnoldi process leads to the following relation
\begin{equation}\label{eqn2021}
(I+\gamma A)^{-1}V_{m}=V_{m}{H}_{m}+{h}_{m+1,m}{\bf v}_{m+1}{\bf e}^{\rm T}_{m},
\end{equation}
where $V_m=[{\bf v}_1,{\bf v}_2,\ldots,{\bf v}_m]$ is an $n\times m$ orthonormal matrix, $n$ is the size of the Toeplitz matrix, ${H}_m=V_m^{\rm T}(I+\gamma A)^{-1}V_m$ is an $m$-by-$m$ upper Hessenberg matrix, and ${\bf e}_{m}$ is the  $m$-th column of the $m$-by-$m$ identity matrix.

Let $\beta=\|{\bf v}\|_2$, if $H_m$ is invertible, then the shift-and-invert Arnoldi method exploits
$$
{\bf y}_{m}(t)=V_{m}\big[{\rm exp}\big(-(t/\gamma)\cdot(H_{m}^{-1}-I)\big)\cdot\beta{\bf e}_{1}\big]\equiv V_{m}{\bf u}_{m}(t)
$$
as an approximation to ${\bf y}(t)$, where ${\bf u}_{m}(t)={\rm
exp}\big(-(t/\gamma)\cdot(H_{m}^{-1}-I)\big)\cdot\beta{\bf e}_{1}$.
The residual is \cite{RRA}
\begin{eqnarray*}\label{eqn206}
{\bf r}_{m}(t)&=&-A{\bf y}_{m}(t)-{\bf y}'_{m}(t)=-AV_{m}{\bf u}_{m}(t)-V_{m}{\bf u}'_{m}(t)\nonumber\\
&=&\frac{h_{m+1,m}}{\gamma}{\bf e}_{m}^{\rm T}H_{m}^{-1}{\bf u}_{m}(t)\cdot(I+\gamma A){\bf v}_{m+1},
\end{eqnarray*}
and
\begin{equation}\label{eqn24}
\|{\bf r}_{m}(t)\|_2=\left|\frac{h_{m+1,m}}{\gamma}{\bf e}_{m}^{\rm
T}H_{m}^{-1}{\bf u}_{m}(t)\right |\cdot\|(I+\gamma A){\bf
v}_{m+1}\|_2,
\end{equation}
which can be used as a cheap stopping criterion in practice.

In the $m$-step shift-and-invert Arnoldi method, we have to compute $m$ Toeplitz matrix-vector products $(I+\gamma A)^{-1}{\bf v}_{i},~i=1,2,\ldots,m$.
Since $\gamma$ is a given shift, we are interested in computing $(I+\gamma A)^{-1}$ once for all. One option is to compute the inverse by some direct methods such as the LU decomposition \cite{GV}.
However, Toeplitz matrix is often dense, and the computation of the inverse of a large dense matrix is prohibitive, especially when the matrix is large.
Fortunately, as $I+\gamma A$ is also a Toeplitz matrix, we have the Gohberg-Semencul formula (GSF) \cite{GSF} for its inverse.
Indeed, the inverse of a Toeplitz matrix $T$ can be reconstructed from its first and last columns. More precisely, denote by ${\bf e}_{1},{\bf e}_{n}$ the first and the last column of the $n$-by-$n$ identity matrix, and let ${\bf x}=[\xi_{0},\xi_{1},\ldots,\xi_{n-1}]^{\rm T}$ and ${\bf y}=[\eta_{0},\eta_{1},\ldots,\eta_{n-1}]^{\rm T}$
be the solutions of the following two Toeplitz systems
\begin{equation}\label{eqn203}
T{\bf x}={\bf e}_{1}\quad {\rm and}\quad  T{\bf y}={\bf e}_{n}.
\end{equation}
If $\xi_{0}\neq0$, then the Gohberg-Semencul formula can be expressed as
\begin{eqnarray}\label{eqn204}
T^{-1}&=&\frac{1}{\xi_{0}}\left\{\left[\begin{array}{cccc}
\xi_{0}&0&\ldots&0\\
\xi_{1}&\xi_{0}&\ldots&0\\
\vdots&\vdots&\ddots&\vdots\\
\xi_{n-1}&\xi_{n-2}&\ldots&\xi_{0}
\end{array}\right]
\left[\begin{array}{cccc}
\eta_{n-1}&\eta_{n-2}&\ldots&\eta_{0}\\
0&\eta_{n-1}&\ldots&\eta_{1}\\
\vdots&\vdots&\ddots&\vdots\\
0&0&\ldots&\eta_{n-1}
\end{array}\right]\right.\nonumber\\
&&-
\left.\left[\begin{array}{cccc}
0&\ldots&0&0\\
\eta_{0}&\ldots&0&0\\
\vdots&\ddots&\vdots&\vdots\\
\eta_{n-2}&\ldots&\eta_{0}&0
\end{array}\right]
\left[\begin{array}{cccc}
0&\xi_{n-1}&\ldots&\xi_{1}\\
\vdots&\vdots&\ddots&\vdots\\
0&0&\ldots&\xi_{n-1}\\
0&0&\ldots&0
\end{array}\right]\right\}\nonumber\\
&\equiv&\frac{1}{\xi_{0}}\left(L_{x}R_{y}-L_{y}^{0}R_{x}^{0}\right),
\end{eqnarray}
where $L_{x},L_{y}^{0}$ are lower Toeplitz matrices, and $R_{y},R_{x}^{0}$ are upper Toeplitz matrices. Consequently, the Toeplitz matrix-vector product $(I+\gamma A)^{-1}{\bf v}_i$ can be realized in several FFTs of length $n$ \cite{LP,Pang}.
We are in a position to present the following algorithm for the Toeplitz matrix exponential; for more details, refer to \cite{LP}.
\begin{algorithm}
{\bf  An shift-and-invert Arnoldi algorithm for product of Toeplitz matrix exponential with a vector}\\
{\bf Step 1.}~Solve the Toeplitz systems $(I+\gamma A){\bf x}={\bf e}_1$ and $(I+\gamma A){\bf y}={\bf e}_n$;\\
{\bf Step 2.}~Choose a convergence tolerance $tol_{\bf exp}$ and the starting vector ${\bf v}_1={\bf v}/\|{\bf v}\|_2$;
\begin{tabbing}
in \= in \= in \= in \= in \= in \= in \= \hspace{3.0truein} \= \kill
\>{\bf for} $i=1,2,\ldots$ do\\
\>\> Perform the shift-and-invert Arnoldi process in which the Toeplitz matrix vector products\\
\>\> $(I+\gamma A)^{-1}{\bf v}_{i}$ are realized through FFTs. \\
\>\>If $\|{\bf r}_i(t)\|_2\leq tol_{\bf exp}$, then form the approximation ${\bf y}_{i}(t)=V_{i}{\bf u}_i(t)$ and Stop, else Continue;\\
\>{\bf end for}
\end{tabbing}
\end{algorithm}

In Step 1 of this algorithm, we have to solve two large scale
non-Hermitian Toeplitz linear systems (\ref{eqn203}). If the desired
accuracy is too high, then we have to pay a large amount of computational cost for solving the Toeplitz linear systems, especially for some ill-conditioned problems.
It is interesting to investigate how to solve the Toeplitz systems
inexactly \cite{LP,Pang}.

\section{An inexact shift-and-invert Arnoldi algorithm for Toeplitz matrix exponential}

\setcounter{equation}{0}

In this section, we consider how to solve the Toeplitz systems
inexactly in the shift-and-invert Arnoldi method. As we solve the Toeplitz linear systems
once for all, it can be understood as an ``inexact" inverse
technology. We first give a new stability analysis on the
Gohberg-Semencul formula with respect to 1-norm and 2-norm, and then establish a relation between the
error of Toeplitz systems and the residual of Toeplitz matrix
exponential. Based on these theoretical results, we propose an inexact shift-and-invert Arnoldi algorithm for Toeplitz matrix exponential.

\subsection{A new stability analysis on the Gohberg-Semencul formula and the GSF condition number}

In this subsection, we give a stability analysis on the Gohberg-Semencul formula and define the ``GSF condition number" of a Toeplitz matrix.
Let $\widetilde{{\bf x}}=[\widetilde{\xi}_{0},\widetilde{\xi}_{1},\ldots,\widetilde{\xi}_{n-1}]^{\rm T}$ and
 $\widetilde{{\bf y}}=[\widetilde{\eta}_{0},\widetilde{\eta}_{1},\ldots,\widetilde{\eta}_{n-1}]^{\rm T}$ be the numerical solutions of $ T{\bf x}={\bf e}_{1}$ and $T {\bf y}={\bf e}_{n}$, respectively. If $\widetilde{\xi}_{0}\neq 0$, we denote
\begin{eqnarray}\label{eqn35}
\widetilde{T}^{-1}&=&\frac{1}{{\widetilde{\xi}}_{0}}\left\{\left[\begin{array}{cccc}
\widetilde{\xi}_{0}&0&\ldots&0\\
\widetilde{\xi}_{1}&\widetilde{\xi}_{0}&\ldots&0\\
\vdots&\vdots&\ddots&\vdots\\
\widetilde{\xi}_{n-1}&\widetilde{\xi}_{n-2}&\ldots&\widetilde{\xi}_{0}
\end{array}\right]
\left[\begin{array}{cccc}
\widetilde{\eta}_{n-1}&\widetilde{\eta}_{n-2}&\ldots&\widetilde{\eta}_{0}\\
0&\widetilde{\eta}_{n-1}&\ldots&\widetilde{\eta}_{1}\\
\vdots&\vdots&\ddots&\vdots\\
0&0&\ldots&\widetilde{\eta}_{n-1}
\end{array}\right]\right.\nonumber\\
&-&
\left.\left[\begin{array}{cccc}
0&\ldots&0&0\\
\widetilde{\eta}_{0}&\ldots&0&0\\
\vdots&\ddots&\vdots&\vdots\\
\widetilde{\eta}_{n-2}&\ldots&\widetilde{\eta}_{0}&0
\end{array}\right]
\left[\begin{array}{cccc}
0&\widetilde{\xi}_{n-1}&\ldots&\widetilde{\xi}_{1}\\
\vdots&\vdots&\ddots&\vdots\\
0&0&\ldots&\widetilde{\xi}_{n-1}\\
0&0&\ldots&0
\end{array}\right]\right\}\nonumber\\
&\equiv&\frac{1}{{\widetilde{\xi}}_{0}}(\widetilde{L}_{x}\widetilde{R}_{y}-\widetilde{L}_{y}^{0}\widetilde{R}_{x}^{0}),
\end{eqnarray}
which is a perturbation to the Toeplitz inverse $T^{-1}$. The following theorem gives an error analysis on the Gohberg-Semencul formula in terms of 1-norm.
\begin{theorem}\label{Thm3.1}
Given $\varepsilon>0$, if $\xi_0\neq 0,~\widetilde{\xi}_0\neq 0$, let $\widetilde{\varepsilon}=\frac{|1/\xi_{0}-1/\widetilde{\xi}_{0}|}{|1/\xi_{0}|}$ be the relative error of $1/\widetilde{\xi}_{0}$ with respect to $1/\xi_0$, and
\begin{equation}\label{3.2}
 \frac{\|\widetilde{{\bf x}}-{\bf x}\|_{1}}{\|{\bf x}\|_{1}}\leq\varepsilon\quad {\rm and}\quad \frac{\|\widetilde{{\bf y}}-{\bf y}\|_{1}}{\|{\bf y}\|_{1}}\leq\varepsilon,
\end{equation}
then we have
 \begin{equation}\label{eqn5.11}
\left\|T^{-1}-\widetilde{T}^{-1}\right\|_{1}\leq\Big|\frac{2}{\xi_{0}}\Big|\cdot\big[\varepsilon+\big(\varepsilon+(1+\varepsilon)
\widetilde{\varepsilon}\big)(1+\varepsilon)\big]\cdot\|{\bf x}\|_{1}\cdot\|{\bf y}\|_{1},
\end{equation}
and
 \begin{equation}\label{eqn3.44}
 \frac{\left\|T^{-1}-\widetilde{T}^{-1}\right\|_{1}}{\|T^{-1}\|_{1}}\leq\Big|\frac{2}{\xi_{0}}\Big|\big[\varepsilon+\big(\varepsilon+(1+\varepsilon)
 \widetilde{\varepsilon}\big)(1+\varepsilon)\big]\cdot{\rm \min}\{\|{\bf x}\|_{1},\|{\bf y}\|_{1}\}.
 \end{equation}
\end{theorem}
{\it Proof.}
 It follows from (\ref{eqn204}) and (\ref{eqn35}) that
\begin{eqnarray}\label{eqn3.6}
\left\|T^{-1}-\widetilde{T}^{-1}\right\|_{1}
&=&\Big\|\frac{1}{\xi_{0}}(L_{x}R_{y}-L_{y}^{0}R_{x}^{0})-\frac{1}{{\widetilde{\xi}_{0}}}(\widetilde{L}_{x}\widetilde{R}_{y}-\widetilde{L}_{y}^{0}\widetilde{R}_{x}^{0})\Big\|_{1}\nonumber\\
&=&\Big\|\Big(\frac{1}{\xi_{0}}L_{x}\Big)R_{y}-L_{y}^{0}\Big(\frac{1}{\xi_{0}}R_{x}^{0}\Big)-\Big(\frac{1}{{\widetilde{\xi}_{0}}}\widetilde{L}_{x}\Big)\widetilde{R}_{y}+\widetilde{L}_{y}^{0}\Big(\frac{1}{{\widetilde{\xi_{0}}}}\widetilde{R}_{x}^{0}\Big)\Big\|_{1}\nonumber\\
&\leq&\Big\|\Big(\frac{1}{\xi_{0}}L_{x}\Big)R_{y}-\Big(\frac{1}{{\widetilde{\xi}_{0}}}\widetilde{L}_{x}\Big)\widetilde{R}_{y}\Big\|_{1}+\Big\|\widetilde{L}_{y}^{0}\Big(\frac{1}{{\widetilde{\xi}_{0}}}\widetilde{R}_{x}^{0}\Big)-L_{y}^{0}\Big(\frac{1}{\xi_{0}}R_{x}^{0}\Big)\Big\|_{1}.
\end{eqnarray}
Moreover, we have that
\begin{eqnarray}\label{eqn3.7}
\Big\|\Big(\frac{1}{\xi_{0}}L_{x}\Big)R_{y}-\Big(\frac{1}{{\widetilde{\xi}_{0}}}\widetilde{L}_{x}\Big)\widetilde{R}_{y}\Big\|_{1}
&=&\Big\|\Big(\frac{1}{\xi_{0}}L_{x}\Big)R_{y}-\Big(\frac{1}{\xi_{0}}L_{x}\Big)\widetilde{R}_{y}+\Big(\frac{1}{\xi_{0}}L_{x}\Big)\widetilde{R}_{y}-\Big(\frac{1}{{\widetilde{\xi}_{0}}}\widetilde{L}_{x}\Big)\widetilde{R}_{y}\Big\|_{1}\nonumber\\
&=&\Big\|\Big(\frac{1}{\xi_{0}}L_{x}\Big)(R_{y}-\widetilde{R}_{y})+\Big(\frac{1}{\xi_{0}}L_{x}-\frac{1}{{\widetilde{\xi}_{0}}}\widetilde{L}_{x}\Big)\widetilde{R}_{y}\Big\|_{1}\nonumber\\
&\leq&\Big|\frac{1}{\xi_{0}}\Big|\cdot\|L_{x}\|_{1}\|R_{y}-\widetilde{R}_{y}\|_{1}+\Big\|\frac{1}{\xi_{0}}L_{x}-\frac{1}{{\widetilde{\xi}_{0}}}\widetilde{L}_{x}\Big\|_{1}\|\widetilde{R}_{y}\|_{1}.
\end{eqnarray}
On the one hand, we obtain from (\ref{3.2}) that
\begin{eqnarray}\label{3110}
\Big\|\frac{1}{\xi_{0}}L_{x}-\frac{1}{{\widetilde{\xi}_{0}}}\widetilde{L}_{x}\Big\|_{1}
&=&\Big\|\frac{1}{\xi_{0}}{\bf x}-\frac{1}{\widetilde{\xi}_{0}}\widetilde{{\bf x}}\Big\|_{1}\nonumber\\
&=&\Big|\frac{1}{\xi_{0}}\Big|\cdot\Big\|{\bf x}-\widetilde{\bf x}+\Big(1-\frac{\xi_{0}}{{\widetilde{\xi}_{0}}}\Big)\widetilde{\bf x})\Big\|_{1}\nonumber\\
&\leq&\Big|\frac{1}{\xi_{0}}\Big|\cdot\big[\|{\bf x}-\widetilde{\bf x}\|_{1}+\widetilde{\varepsilon}\|\widetilde{\bf x}\|_{1}\big]\nonumber\\
&\leq&\Big|\frac{1}{\xi_{0}}\Big|\cdot\big[\varepsilon+(1+\varepsilon)\widetilde{\varepsilon}\big]\cdot\|{\bf x}\|_{1},
\end{eqnarray}
where $\widetilde{\varepsilon}=\frac{|1/\xi_{0}-1/\widetilde{\xi}_{0}|}{|1/\xi_{0}|}$ is the relative error of $1/\xi_{0}$.
On the other hand, we note from (\ref{eqn204}) and (\ref{eqn35}) that
\begin{equation}\label{eqn3.8}
\|L_{x}\|_{1}=\|{\bf x}\|_{1},~\|\widetilde{R}_{y}\|_{1}=\|\widetilde{\bf y}\|_{1}\leq(1+\varepsilon)\|{\bf y}\|_{1},
\end{equation}
and
\begin{equation}\label{399}
\|R_{y}-\widetilde{R}_{y}\|_{1}=\|{\bf y}-\widetilde{\bf y}\|_{1}\leq\varepsilon\|{\bf y}\|_{1}.
\end{equation}
From (\ref{eqn3.7})--(\ref{399}), we obtain
\begin{eqnarray}\label{eqn3.11}
\Big\|\Big(\frac{1}{\xi_{0}}L_{x}\Big)R_{y}-\Big(\frac{1}{\widetilde{\xi_{0}}}\widetilde{L}_{x}\Big)\widetilde{R}_{y}\Big\|_{1}
&\leq&\Big|\frac{1}{\xi_{0}}\Big|\cdot\|{\bf x}\|_{1}\cdot\|{\bf y}\|_{1}\cdot\varepsilon+\Big|\frac{1}{\xi_{0}}\Big|\cdot\big[\varepsilon+(1+\varepsilon)\widetilde{\varepsilon}\big]\cdot\|{\bf x}\|_{1}\cdot\|{\bf y}\|_{1}\cdot(1+\varepsilon)\nonumber\\
&=&\Big|\frac{1}{\xi_{0}}\Big|\cdot\big[\varepsilon+\big(\varepsilon+(1+\varepsilon)\widetilde{\varepsilon}\big)(1+\varepsilon)\big]\cdot\|{\bf x}\|_{1}\cdot\|{\bf y}\|_{1}.
\end{eqnarray}
Similarly, for the second part of (\ref{eqn3.6}), we can prove that
\begin{equation}\label{eqn5.10}
\Big\|\widetilde{L}_{y}^{0}\Big(\frac{1}{\widetilde{\xi_{0}}}\widetilde{R}_{x}^{0}\Big)-L_{y}^{0}\Big(\frac{1}{\xi_{0}}R_{x}^{0}\Big)\Big\|_{1}\leq\Big|\frac{1}
{\xi_{0}}\Big|\cdot\big[\varepsilon+\big(\varepsilon+(1+\varepsilon)\widetilde{\varepsilon}\big)(1+\varepsilon)\big]\cdot\|{\bf x}\|_{1}\cdot\|{\bf y}\|_{1}.
\end{equation}
Combining (\ref{eqn3.6}), (\ref{eqn3.11}) and (\ref{eqn5.10}), we arrive at
$$
\left\|T^{-1}-\widetilde{T}^{-1}\right\|_{1}\leq\Big|\frac{2}{\xi_{0}}\Big|\cdot\big[\varepsilon+\big(\varepsilon+(1+\varepsilon)
\widetilde{\varepsilon}\big)(1+\varepsilon)\big]\cdot\|{\bf x}\|_{1}\cdot\|{\bf y}\|_{1}.
$$
By $(\ref{eqn203})$, we have that
 $ \|{\bf x}\|_{1}\leq\|T^{-1}\|_1$ and $\|{\bf y}\|_{1}\leq\|T^{-1}\|_{1}$.
Thus,
$$
\left\|T^{-1}-\widetilde{T}^{-1}\right\|_1\leq\Big|\frac{2}{\xi_{0}}\Big|\cdot\big[\varepsilon+\big(\varepsilon+(1+\varepsilon)
\widetilde{\varepsilon}\big)(1+\varepsilon)\big]\cdot\|{\bf x}\|_1\cdot\|T^{-1}\|_1,
 $$
and
$$
 \left\|T^{-1}-\widetilde{T}^{-1}\right\|_1\leq\Big|\frac{2}{\xi_{0}}\Big|\cdot\big[\varepsilon+\big(\varepsilon+
 (1+\varepsilon)\widetilde{\varepsilon}\big)(1+\varepsilon)\big]\cdot\|{\bf y}\|_1\cdot\|T^{-1}\|_1,
$$
a combination of which yields (\ref{eqn3.44}).
$\hfill \Box$


Furthermore, we have the following corollary on the relative error of Toeplitz inverse.

\begin{Cor}\label{Col3.1}
Under the above notations, there holds
\begin{equation}\label{eqn3146}
 \frac{\left\|T^{-1}-\widetilde{T}^{-1}\right\|_{1}}{\|T^{-1}\|_{1}}\leq\frac{\|T\|_1\|{\bf y}\|_1}{|\xi_{0}|/\|{\bf x}\|_1}\cdot 2\big[\varepsilon+\big(\varepsilon+(1+\varepsilon)
 \widetilde{\varepsilon}\big)(1+\varepsilon)\big].
\end{equation}
\end{Cor}

{\it Proof.}
We note that $\|T\|_{1}\|{\bf y}\|_{1}\geq \|T{\bf y}\|_{1}= 1$ and $\|T\|_{1}\|{\bf x}\|_{1}\geq 1$. It follows from Theorem \ref{Thm3.1} that
\begin{eqnarray*}
 \frac{\left\|T^{-1}-\widetilde{T}^{-1}\right\|_{1}}{\|T^{-1}\|_{1}}&\leq&\Big|\frac{2}{\xi_{0}}\Big|\big[\varepsilon+\big(\varepsilon+(1+\varepsilon)
 \widetilde{\varepsilon}\big)(1+\varepsilon)\big]\cdot{\rm \min}\{\|{\bf x}\|_{1},\|{\bf y}\|_{1}\}\nonumber\\
 &\leq&\frac{\|T\|_{1}\|{\bf y}\|_{1}\|{\bf x}\|_{1}}{|\xi_{0}|}\cdot 2\big[\varepsilon+\big(\varepsilon+(1+\varepsilon)\widetilde{\varepsilon}\big)(1+\varepsilon)\big]\nonumber\\
&=&\frac{\|T\|_{1}\|{\bf y}\|_{1}}{|\xi_{0}|/\|{\bf x}\|_{1}}\cdot 2\big[\varepsilon+\big(\varepsilon+(1+\varepsilon)\widetilde{\varepsilon}\big)(1+\varepsilon)\big].
 \end{eqnarray*}
$\hfill \Box$

By {\rm(}\ref{eqn3146}{\rm)}, $(\|T\|_1\|{\bf y}\|_1)/(|\xi_{0}|/\|{\bf x}\|_{1})$ is an enlarge factor of the solution of $T^{-1}$ over the vector error $\varepsilon$. So we can give the following definition on the condition number of a Toeplitz matrix.
\begin{definition}\label{Def3.1}
We define
\begin{equation}\label{eqn3155}
\kappa_1^{GSF}(T)=\frac{\|T\|_1\|{\bf y}\|_1}{|\xi_{0}|/\|{\bf x}\|_1}
\end{equation}
as the 1-norm ``GSF condition number" of a Toeplitz matrix.
\end{definition}

Note that $|\xi_0|/\|{\bf x}\|_1$ is the ``proportion" of $|\xi_0|$ with respect to $\|{\bf x}\|_1$, and
$$
\kappa_{\rm eff}^{\rm I}(T)=\frac{\|T\|_1\|{\bf y}\|_1}{\|{\bf e}_n\|_1}=\|T\|_1\|{\bf y}\|_1
$$
is the ``effective" 1-norm condition number of $T{\bf y}={\bf e}_n$ defined in \cite{CMSW}.
Moreover, we notice that
$$
\kappa_{\rm eff}^{\rm II}(T)=\frac{\left\|T^{-1}\right\|_1\|{\bf e}_n\|_1}{\|{\bf y}\|_1}
$$
the ``effective" 1-norm condition number of $T{\bf y}={\bf e}_n$ defined in \cite{Li,Rice}, and we have that
\begin{equation}\label{eqn316}
\kappa_1(T)=\kappa_{\rm eff}^{\rm I}(T)\cdot\kappa_{\rm eff}^{\rm II}(T),
\end{equation}
where $\kappa_1(T)=\|T\|_1\left\|T^{-1}\right\|_1$ is the ``classical" 1-norm condition number \cite{GV} of the matrix $T$.

\begin{rem}
In terms of Corollary \ref{Col3.1}, $\kappa_1^{GSF}(T)$ is an estimation to $\kappa_1(T)$. By {\rm(}\ref{eqn316}{\rm)}, $\|{\bf x}\|_1/|\xi_0|$ can be used as an approximation to $\kappa_{\rm eff}^{\rm II}(T)$.
Thus, an advantage of {\rm(}\ref{eqn3155}{\rm)} is that one can evaluate the ``classical" condition number $\kappa_1(T)$ of a Toeplitz matrix, and the ``effective" condition numbers $\kappa_{\rm eff}^{\rm I}(T),\kappa_{\rm eff}^{\rm II}(T)$ via solving Toeplitz systems, with no need to form the Toeplitz inverse explicitly.
\end{rem}

Since $V_m$ is orthonormal, it is desirable to investigate the absolute and relative errors of $\widetilde{T}^{-1}$ with respect to $T^{-1}$ according to 2-norm. We have the following result.
\begin{theorem}\label{the5.1}
Under the assumptions and notations of Theorem \ref{Thm3.1}, we have that
 \begin{equation}\label{eqn3.4}
\left \|T^{-1}-\widetilde{T}^{-1}\right
\|_{2}\leq\left|\frac{2}{\xi_{0}}\right|\big[\varepsilon+\big(\varepsilon+(1+\varepsilon)
 \widetilde{\varepsilon}\big)(1+\varepsilon)\big]\cdot\|{\bf x}\|_{1}\|{\bf y}\|_{1},
 \end{equation}
 and
 \begin{equation}\label{eqn325}
 \frac{\left\|T^{-1}-\widetilde{T}^{-1}\right\|_{2}}{\|T^{-1}\|_{2}}\leq\left|\frac{2\sqrt{n}}{\xi_{0}}\right|\big[\varepsilon+\big(\varepsilon+(1+\varepsilon)
 \widetilde{\varepsilon}\big)(1+\varepsilon)\big]\cdot{\rm \min}\{\|{\bf x}\|_{1},\|{\bf y}\|_{1}\}.
 \end{equation}
\end{theorem}

{\it Proof.}
Recall that
\begin{equation}\label{39}
\left\|T^{-1}-\widetilde{T}^{-1}\right\|_{2}^{2}\leq
\left\|T^{-1}-\widetilde{T}^{-1}\right\|_{1}\cdot
\left\|T^{-1}-\widetilde{T}^{-1}\right\|_{\infty}.
\end{equation}
On the one hand, we have from (\ref{eqn5.11}) that
$$
\left\|T^{-1}-\widetilde{T}^{-1}\right\|_{1}\leq\left|\frac{2}{\xi_{0}}\right|\cdot\big[\varepsilon+(1+\varepsilon)\big(\varepsilon+(1+\varepsilon)
\widetilde{\varepsilon}\big)\big]\cdot\|{\bf x}\|_{1}\|{\bf
y}\|_{1}.
$$
On the other hand, we can give an upper bound on $\left\|T^{-1}-\widetilde{T}^{-1}\right\|_{\infty}$:
\begin{equation}\label{eqn333}
\left\|T^{-1}-\widetilde{T}^{-1}\right\|_{\infty}\leq\Big|\frac{2}{\xi_{0}}\Big|\cdot\big[\varepsilon+\big(\varepsilon+(1+\varepsilon)
\widetilde{\varepsilon}\big)(1+\varepsilon)\big]\cdot\|{\bf x}\|_{1}\|{\bf y}\|_{1},
\end{equation}
whose proof is similar to that of Theorem 3.1, see \cite{Our}.
So we have from (\ref{39}), (\ref{eqn5.11}) and
(\ref{eqn333}) that
$$
\left\|T^{-1}-\widetilde{T}^{-1}\right\|_{2}\leq\left|\frac{2}{\xi_{0}}\right|\cdot\big[\varepsilon+(1+\varepsilon)\big(\varepsilon+(1+\varepsilon)
\widetilde{\varepsilon}\big)\big]\cdot\|{\bf x}\|_{1}\|{\bf
y}\|_{1}.
$$

For (\ref{eqn325}), we notice that
\begin{equation}\label{320}
\|{\bf x}\|_{1}\leq \sqrt{n}~\|{\bf x}\|_{2}= \sqrt{n}~
\left\|T^{-1} {\bf  e}_1 \right\|_{2} \leq \sqrt{n}~
\left\|T^{-1}\right\|_{2}\|{\bf
e}_1\|_{2}=\sqrt{n}~\left\|T^{-1}\right\|_{2}.
\end{equation}
Combining (\ref{eqn3.4}) and (\ref{320}), we drive
$$
 \frac{\left\|T^{-1}-\widetilde{T}^{-1}\right\|_{2}}{\left\|T^{-1}\right\|_{2}}\leq\left|\frac{2}{\xi_{0}}\right|\big[\varepsilon+(1+\varepsilon)\big(\varepsilon+(1+\varepsilon)
 \widetilde{\varepsilon}\big)\big]\cdot\sqrt{n}~\|{\bf y}\|_{1}.
$$
Similarly, we can prove that
$$
 \frac{\left\|T^{-1}-\widetilde{T}^{-1}\right\|_{2}}{\left\|T^{-1}\right\|_{2}}\leq\left|\frac{2}{\xi_{0}}\right|\big[\varepsilon+(1+\varepsilon)\big(\varepsilon+(1+\varepsilon)
 \widetilde{\varepsilon}\big)\big]\cdot\sqrt{n}~\|{\bf x}\|_{1},
$$
a combination of the above two inequalities yields (\ref{eqn325}).
$\hfill \Box$

\begin{rem}
In {\rm(}\cite[p.321]{MHG}{\rm)}, Gutknecht and Hochbruck analyzed
the stability of the Gohberg-Semencul formula, and gave the
following two upper bounds for the absolute and relative errors with respect to
$T^{-1}$:
\begin{equation}\label{eqn3.14}
\left\|T^{-1}-\widetilde{T}^{-1}\right\|_{2}\leq\left|\frac{1}{\xi_{0}}\right|\Big(\|{\bf
x}\|_{2}\|{\bf y}\|_{2}\left(4n \varepsilon+2n^{2}
\epsilon+2n\epsilon\right)+\sqrt{n}\epsilon
\left\|T^{-1}\right\|_{2}\Big),
 \end{equation}
 and
 \begin{equation}\label{eqn3156}
\frac{\left\|T^{-1}-\widetilde{T}^{-1}\right\|_{2}}{\left\|T^{-1}\right\|_{2}}
\leq
\left|\frac{1}{\xi_{0}}\right|\Big(2n\left\|T^{-1}\right\|_{2}(2\varepsilon+(n+1)\epsilon\big)+\sqrt{n}\epsilon\Big),
\end{equation}
where $n$ is the size of the Toeplitz matrix, $\epsilon$ is the machine precision and $\varepsilon$ is a normwise relative error bound
that satisfies
$$
\left\|\widetilde{{\bf x}}\right\|_{2}\leq \|{\bf
x}\|_2\cdot(1+\varepsilon),\quad \left\|\widetilde{{\bf
y}}\right\|_{2}\leq \|{\bf y}\|_2\cdot(1+\varepsilon).
$$
By setting $\epsilon=0$, {\rm(}\ref{eqn3.14}{\rm)} and {\rm(}\ref{eqn3156}{\rm)} reduce to
\begin{equation}\label{eqn.3.13}
\left\|T^{-1}-\widetilde{T}^{-1}\right\|_{2} \leq
\left|\frac{4n}{\xi_{0}}\right|\cdot\|{\bf x}\|_{2}\|{\bf
y}\|_{2}\cdot \varepsilon,
\end{equation}
and
\begin{equation}\label{eqn324}
\frac{\left\|T^{-1}-\widetilde{T}^{-1}\right\|_{2}}{\left\|T^{-1}\right\|_{2}}
\leq\left|\frac{4
n}{\xi_{0}}\right|\cdot\left\|T^{-1}\right\|_{2}\cdot\varepsilon,
\end{equation}
respectively.

Compared with {\rm(}\ref{eqn.3.13}{\rm)} and
{\rm(}\ref{eqn324}{\rm)},
the new bounds given in {\rm(}\ref{eqn5.11}{\rm)} and {\rm(}\ref{eqn3.4}{\rm)} are independent of $n$,
while {\rm(}\ref{eqn325}{\rm)} relies on
$\sqrt{n}$ instead of $n$. Thus, our new upper bounds can be  sharper
than {\rm(}\ref{eqn.3.13}{\rm)} and {\rm(}\ref{eqn324}{\rm)} when
$n$ is large and $\widetilde{\varepsilon}=\mathcal{O}(\varepsilon)$.
See Example 3 of Section 4 for a comparison of these upper bounds.
\end{rem}

\subsection{Relationship between the error of Toeplitz
systems and the residual of Toeplitz matrix exponential}

In this subsection, we establish a relationship between the error of Toeplitz
systems and the residual of Toeplitz matrix exponential, and propose an inexact shift-and-invert Arnoldi method for product of a Toeplitz matrix exponential with a vector. It is shown
that if the GSF condition number of the Toeplitz matrix is medium
sized, we can solve the Toeplitz systems in a relatively low accuracy.

For
simplicity, in the following we denote $T=I+\gamma A$ whenever
necessary. Indeed, if the Toeplitz systems (\ref{eqn203}) are solved
inexactly, the errors of the matrix-vector products can be expressed
as ${\bf f}_{i}=\widetilde{T}^{-1}{\bf
v}_{i}-T^{-1}{\bf v}_{i},~i=1,2,\ldots,m$. Let $F_{m}=[{\bf f}_{1},{\bf
f}_{2},\ldots,{\bf f}_{m}]$, we get the following relation
for the $m$-step ``inexact" shift-and-invert Arnoldi procedure
\begin{equation}\label{eqn308}
(I+\gamma A)^{-1}V_{m}+F_{m}=V_{m}\widetilde{H}_{m}+\widetilde{h}_{m+1,m}{\bf v}_{m+1}{\bf e}_{m}^{\rm T},
\end{equation}
where $V_{m}=[{\bf v}_{1},{\bf v}_{2},\ldots,{\bf v}_{m}]$ is an $n\times m$ orthonormal matrix, and $\widetilde{H}_{m}=V_{m}^{\rm T}\big[(I+\gamma A)^{-1}+F_mV_m^{\rm T}\big]V_{m}$ is an $m\times m$ upper Hessenberg matrix. Note that $V_m$ is different from the one given in (\ref{eqn2021}), and the subspace spanned by $V_m$ is not a Krylov subspace any more.
\begin{lemma}\label{the5.1}
If $\widetilde{H}_{m}$ is invertible, denote
$G=-\frac{1}{\gamma}(I+\gamma
A)F_{m}\widetilde{H}_{m}^{-1}V_{m}^{\rm T}$, then the ``inexact"
shift-and-invert Arnoldi relation {\rm(}\ref{eqn308}{\rm)} can be
rewritten as
\begin{equation}
\big(A+G\big)V_{m}=V_{m}\left(\frac{1}{\gamma}\left(\widetilde{H}_{m}^{-1}-I\right)\right)-\frac{\widetilde{h}_{m+1,m}}{\gamma}(I+\gamma
A){\bf v}_{m+1}{\bf e}_{m}^{\rm T}\widetilde{H}_{m}^{-1}.
\end{equation}
\end{lemma}

{\it Proof.}
Multiplying $I+\gamma A$ on both sides of (\ref{eqn308}) yields
$$
V_{m}+(I+\gamma A)F_{m}=(I+\gamma A)V_{m}\widetilde{H}_{m}+\widetilde{h}_{m+1,m}(I+\gamma A){\bf v}_{m+1}{\bf e}_{m}^{\rm T},
$$
that is,
\begin{eqnarray*}
AV_{m}-\frac{1}{\gamma}(I+\gamma
A)F_{m}\widetilde{H}_{m}^{-1}
&=&\frac{1}{\gamma}V_{m}\left(I-\widetilde{H}_{m}\right)\widetilde{H}_{m}^{-1}-\frac{\widetilde{h}_{m+1,m}}{\gamma}(I+\gamma
A){\bf v}_{m+1}{\bf e}_{m}^{\rm T}\widetilde{H}_{m}^{-1}\\\nonumber
&=&\frac{1}{\gamma}V_{m}\left(\widetilde{H}_{m}^{-1}-I\right)-\frac{\widetilde{h}_{m+1,m}}{\gamma}(I+\gamma
A){\bf v}_{m+1}{\bf e}_{m}^{\rm T}\widetilde{H}_{m}^{-1}\\\nonumber
&=&V_{m}\left(\frac{1}{\gamma}\left(\widetilde{H}_{m}^{-1}-I\right)\right)-\frac{\widetilde{h}_{m+1,m}}{\gamma}(I+\gamma
A){\bf v}_{m+1}{\bf e}_{m}^{\rm T}\widetilde{H}_{m}^{-1}.\nonumber
\end{eqnarray*}
The above equation can be rewritten as
$$
\left(A-\frac{1}{\gamma}(I+\gamma
A)F_{m}\widetilde{H}_{m}^{-1}V_{m}^{\rm
T}\right)V_{m}=\big(A+G\big)V_{m}=V_{m}\left(\frac{1}{\gamma}\left(\widetilde{H}_{m}^{-1}-I\right)\right)-\frac{\widetilde{h}_{m+1,m}}{\gamma}(I+\gamma
A){\bf v}_{m+1}{\bf e}_{m}^{\rm T}\widetilde{H}_{m}^{-1}.
$$
$\hfill \Box$

Let ${\bf u}_{m}(t)={\rm exp}\big(-(t/\gamma)\cdot(\widetilde{H}_{m}^{-1}-I)\big)\cdot\beta{\bf e}_{1}$, then we can use
${\bf y}_m(t)=V_m{\bf u}_{m}(t)$ as an approximation to ${\bf y}(t)$. The ``real" residual is defined as \cite{RRA}
\begin{equation}\label{eqn320}
{\bf r}^{real}=-AV_{m}{\bf u}_{m}(t)-V_{m}{\bf u}'_{m}(t).
\end{equation}
However, it is not computable since $AV_m$ is unavailable in
practice. Thus, we define
\begin{eqnarray}\label{eqn310}
{\bf r}^{comp}&=&-(A+G)V_{m}{\bf u}_{m}(t)-V_{m}{\bf
u}'_{m}(t)\\\nonumber
&=&-V_{m} \frac{\left(\widetilde{H}_{m}^{-1}-I\right)}{\gamma}{\bf u}_{m}(t)+\frac{\widetilde{h}_{m+1,m}}{\gamma}(I+\gamma
A){\bf v}_{m+1}{\bf e}_{m}^{\rm T}\widetilde{H}_{m}^{-1}{\bf
u}_{m}(t)+V_{m} \frac{\left(\widetilde{H}_{m}^{-1}-I\right)}{\gamma}{\bf
u}_{m}(t)\\\nonumber
&=&\frac{\widetilde{h}_{m+1,m}}{\gamma}(I+\gamma A){\bf v}_{m+1}{\bf
e}_{m}^{\rm T}\widetilde{H}_{m}^{-1}{\bf u}_{m}(t),
\end{eqnarray}
as the ``computed" residual.
Moreover,
\begin{equation}\label{eqn330}
\|{\bf
r}^{comp}\|_{2}=\left|\frac{\widetilde{h}_{m+1,m}}{\gamma}{\bf
e}^{\rm T}_{m}\widetilde{H}_{m}^{-1}{\bf
u}_{m}(t)\right|\cdot\|(I+\gamma A){\bf v}_{m+1}\|_{2},
\end{equation}
which can be used as a cheap stopping criterion in the ``inexact" shift-and-invert Arnoldi method for Toeplitz matrix exponential.

We are ready to provide a practical stopping criterion for solving
the Toeplitz systems inexactly. The key is how to investigate the
distance between ${\bf r}^{real}$ and ${\bf r}^{comp}$. It is seen
that
\begin{eqnarray}\label{eqn311}
\left\|{\bf r}^{real}-{\bf r}^{comp}\right\|_{2}&=&\|GV_{m}{\bf u}_{m}(t)\|_{2}\nonumber\\
&=&\left\|\frac{1}{\gamma}(I+\gamma A)F_{m}\widetilde{H}_{m}^{-1}{\bf u}_{m}(t)\right\|_{2}\nonumber\\
&\leq&\left|\frac{1}{\gamma}\right|\cdot\|I+\gamma A\|_{2}\cdot
\left\|F_{m}\widetilde{H}_{m}^{-1}{\bf u}_{m}(t)\right\|_{2}.
\end{eqnarray}

From (\ref{eqn3.4}), we obtain
\begin{eqnarray*}\label{eqn313}
\|{\bf f}_{i}\|_{2}&=&\left\|\widetilde{T}^{-1}{\bf
v}_{i}-T^{-1}{\bf v}_{i}\right\|_{2}\leq
\left\|T^{-1}-\widetilde{T}^{-1}\right\|_{2}\|{\bf
v}_{i}\|_{2}\\\nonumber
&=&\left\|T^{-1}-\widetilde{T}^{-1}\right\|_{2} \leq
\left|\frac{2}{\xi_{0}}\right|\cdot\big[\varepsilon+\big(\varepsilon+(1+\varepsilon)\widetilde{\varepsilon}\big)(1+\varepsilon)\big]\cdot\|{\bf
x}\|_{1} \|{\bf y}\|_{1},\qquad i=1,2,\ldots,m.
\end{eqnarray*}
If $\widetilde{\varepsilon}\leq\varepsilon\ll1$,
then $\varepsilon+(\varepsilon+(1+\varepsilon)\widetilde{\varepsilon})(1+\varepsilon)\leq 3\varepsilon+\mathcal{O}(\varepsilon^{2})$, and
\begin{eqnarray}\label{eqn313}
\|{\bf f}_{i}\|_{2}&\leq
&\left|\frac{2}{\xi_{0}}\right|\cdot\big[\varepsilon+\big(\varepsilon+(1+\varepsilon)\widetilde{\varepsilon}\big)(1+\varepsilon)\big]\cdot\|{\bf
x}\|_{1}\|{\bf y}\|_{1}\\\nonumber
&\lesssim&\left|\frac{6}{\xi_{0}}\right|\cdot\|{\bf x}\|_{1}\|{\bf
y}\|_{1}\cdot\varepsilon,\qquad i=1,2,\ldots,m,
\end{eqnarray}
where we omit the high order term $\mathcal{O}(\varepsilon^2)$.

Denote $\widetilde{H}_{m}^{-1}{\bf
u}_{m}(t)=[\alpha_{1},\alpha_{2},\ldots,\alpha_{m}]^{\rm T}$, then
\begin{eqnarray}\label{eqn3141}
\left\|F_{m}\widetilde{H}_{m}^{-1}{\bf
u}_{m}(t)\right\|_{2}&=&\left\|\sum_{i=1}^{m} {\bf
f}_{i}\alpha_{i}\right\|_2\leq \max _{1\leq i\leq m}\|{\bf
f}_{i}\|_2\cdot \left\|\widetilde{H}_{m}^{-1}{\bf
u}_{m}(t)\right\|_{1}\\\nonumber &\leq& \max _{1\leq i\leq m}\|{\bf
f}_{i}\|_2\cdot\sqrt{m}~\left\|\widetilde{H}_{m}^{-1}{\bf
u}_{m}(t)\right\|_{2}.
\end{eqnarray}
Combining (\ref{eqn311}), (\ref{eqn313}) and (\ref{eqn3141}), we have that
\begin{eqnarray*}\label{eqn315}
\left\|{\bf r}^{real}-{\bf
r}^{comp}\right\|_{2}&\leq&\left|\frac{1}{\gamma}\right|\|I+\gamma
A\|_{2}\cdot\max _{1\leq i\leq m}\|{\bf
f}_{i}\|_2\cdot\sqrt{m}~\left\|\widetilde{H}_{m}^{-1}\right\|_{2}\|{\bf
u}_{m}(t)\|_{2}\\\nonumber
&\lesssim&\left|\frac{1}{\gamma}\right|\|I+\gamma A\|_{2}\cdot
\left|\frac{6}{\xi_{0}}\right|\|{\bf x}\|_{1}\|{\bf
y}\|_{1}\varepsilon\cdot\sqrt{m}~\left\|\widetilde{H}_{m}^{-1}\right\|_{2}\|{\bf
u}_{m}(t)\|_{2}.\nonumber
\end{eqnarray*}

Let $tol_{\bf exp}$ be the convergence threshold for the
shift-and-invert Arnoldi method for solving (\ref{eqn11}). If
$$
\varepsilon\leq\frac{|\gamma\xi_{0}|\cdot tol_{\bf
exp}}{6\sqrt{m}\|I+\gamma A\|_{2}\cdot\|{\bf x}\|_{1}\|{\bf
y}\|_{1}\cdot \left\|\widetilde{H}_{m}^{-1}\right\|_{2}\|{\bf
u}_{m}(t)\|_{2}},
$$
then we have that
$$
\|{\bf r}^{real}-{\bf r}^{comp}\|_{2} \lesssim tol_{\bf exp}.
$$
In conclusion, we
have the following theorem.
\begin{theorem}\label{Thm3.3}
Under the above notations and assumptions, if
\begin{eqnarray}
\varepsilon &\leq& \frac{|\xi_{0}|/\|{\bf x}\|_{1}}{\|I+\gamma A\|_{2}\|{\bf y}\|_{1}}\cdot\frac{|\gamma|\cdot tol_{\bf exp}}{6\sqrt{m}\cdot
\left\|\widetilde{H}_{m}^{-1}\right\|_{2}\|{\bf u}_{m}(t)\|_{2}}\nonumber\\
&=&\frac{|\xi_{0}|/\|{\bf x}\|_{1}}{\|I+\gamma A\|_{1}\|{\bf
y}\|_{1}}\cdot\frac{|\gamma|\cdot\|I+\gamma A\|_{1}\cdot tol_{\bf
exp}}{6\sqrt{m}~ \|I+\gamma A\|_{2}\cdot
\left\|\widetilde{H}_{m}^{-1}\right\|_{2}\|{\bf
u}_{m}(t)\|_{2}}\label{eqn343}
\end{eqnarray}
then
$$
\|{\bf r}^{real}-{\bf r}^{comp}\|_{2}\lesssim tol_{\bf exp}.
$$
\end{theorem}

\begin{rem}
We notice that
$$
\kappa_1^{GSF}(I+\gamma A)=\frac{\|I+\gamma A\|_{1}\|{\bf y}\|_{1}}{|\xi_{0}|/\|{\bf x}\|_{1}}
$$
is just the 1-norm ``GSF condition number" defined in Definition \ref{Def3.1}, which can be utilized as an estimation to the 1-norm condition number of $I+\gamma A$.
Furthermore, {\rm(}\ref{eqn343}{\rm)} can be reformulated as
\begin{equation}
\varepsilon \leq\frac{1}{\kappa_1^{GSF}(I+\gamma
A)}\cdot\frac{|\gamma|\cdot\|I+\gamma A\|_{1}\cdot tol_{\bf
exp}}{6\sqrt{m}~\|I+\gamma A\|_{2} \cdot
\left\|\widetilde{H}_{m}^{-1}\right\|_{2}\|{\bf
u}_{m}(t)\|_{2}}\label{eqn335}.
\end{equation}
This implies that if the GSF condition number of the Toeplitz matrix is medium sized, we can solve the Toeplitz systems in a {\rm(}relatively{\rm)} low accuracy. Otherwise, we have to
solve the Toeplitz systems in a {\rm(}relatively{\rm)} high accuracy.
\end{rem}

\begin{rem}
Unfortunately, the parameters $\left\|\widetilde{H}_{m}^{-1}\right\|_{2}$,
$\|{\bf u}_{m}(t)\|_{2}$, and $\kappa_1^{GSF}(I+\gamma A)$ are
unavailable {\it a prior}. We notice that
$\left\|\widetilde{H}_{m}^{-1}\right\|_{2}$ is uniformly bounded and
$\|{\bf u}_{m}(t)\|_{2}=\mathcal{O}(\|{\bf y}(t)\|_2)$ as the
shift-and-invert Arnoldi method converges. Therefore, if
$\kappa_1^{GSF}(I+\gamma A)$ is medium sized and $\|I+\gamma
A\|_{1}=\mathcal{O}\left(\left\|\widetilde{H}_{m}^{-1}\right\|_{2}\|{\bf
u}_{m}(t)\|_{2}\right)$, we suggest using
 \begin{equation}\label{eqn5.21}
\left\{\|{\bf r}_x\|_2,\|{\bf r}_y\|_2\right\}
\leq\frac{|\gamma|}{6\sqrt{m}\cdot\max\{\|{\bf fcol}\|_{2},\|{\bf
frow}\|_{2}\}}\cdot{tol_{\bf exp}}
 \end{equation}
as the stopping criterion for solving the Toeplitz systems, where ${\bf r}_x={\bf e}_1-T\widetilde{{\bf x}}$ and ${\bf r}_y={\bf e}_n-T\widetilde{{\bf y}}$
are the residuals of the Toeplitz systems, and {\bf fcol}, {\bf frow} are the first column and first row of $I+ \gamma A$, respectively.
\end{rem}

In summary, we propose the following ``inexact" shift-and-invert Arnoldi algorithm for solving the Toeplitz matrix exponential problem (\ref{eqn11}).
\begin{algorithm}
{\bf ~An inexact shift-and-invert Arnoldi algorithm for product of Toeplitz matrix exponential with a vector}\\
 This algorithm is similar to Algorithm 1, except for the Toeplitz linear systems {\rm(}\ref{eqn203}{\rm)} are solved inexactly in Step 1,
with the stopping criterion described in {\rm(}\ref{eqn5.21}{\rm)}.
\end{algorithm}

We point out that several results on computing matrix functions using inexact Krlov methods have been developed in other contexts, where 2-norm estimates are given.
In \cite{Frommer}, Frommer {\it et al.} considered how to cheaply recover a secondary Lanczos process starting at an arbitrary
Lanczos vector. This secondary process is then used to efficiently obtain computable error
estimates and error bounds for the Lanczos approximations to the action of a rational matrix function
on a vector, e.g., the matrix sign function.

\section{Numerical experiments}

\setcounter{equation}{0}


In this section, we perform some numerical examples to show the
efficiency of Algorithm 2 and the effectiveness of our theoretical
results. All the numerical experiments were run on two core Intel(R)
Core(TM)2 E7400 processor with CPU 2.8 GHz and RAM 1.99 GB, under
the Windows 7 operating system. The experimental results were
obtained by using a MATLAB 7.7 implementation with machine precision
$\epsilon\approx 2.22\times 10^{-16}$.

two core Intel(R) Core(TM)2 E7400 processor with CPU 2.8 GHz and RAM
2 GB

As was done in \cite{LP}, we use the (unrestarted) GMRES algorithm
\cite{Saad2} with T. Chan's optimal (circulant) preconditioner \cite{CJ,TChan,Ng2} for solving
the Toeplitz systems in Algorithm 1 and Algorithm 2.  Let $tol_{\bf exp}$, $tol_{\bf sys}$ be the tolerance for computing the Toeplitz matrix exponential-vector product, and that for solving the Toeplitz systems, respectively. Denote by
$\mathcal{M}$ the optimal preconditioner due to T. Chan, by $\widetilde{\bf
q}=\widetilde{\bf x}$ (or $\widetilde{\bf y}$) the approximate
solution of the Toeplitz system, and by ${\bf b}={\bf e}_{1}$ (or ${\bf e}_{n}$) the
right-hand side. In Algorithm 2 we use
\begin{equation}\label{eqn42}
\left\|\mathcal{M}^{-1}{\bf b}-\mathcal{M}^{-1}(I+\gamma
A)\widetilde{\bf q}\right\|_2\leq\frac{|\gamma|}{6\sqrt{100}\cdot
\max\{\|{\bf fcol}\|_{2},\|{\bf frow}\|_{2}\}}\cdot tol_{\bf
exp}\equiv tol_{\bf sys}
\end{equation}
as the stopping criterion for the Toeplitz systems, where {\bf fcol} and {\bf frow} denote the first column and the first row of $I+ \gamma A$, respectively. This algorithm mimics solving the two Toeplitz systems ``inexactly" with an iterative solver.

In Algorithm 1, we use
\begin{equation}
\left\|\mathcal{M}^{-1}{\bf b}-\mathcal{M}^{-1}(I+\gamma
A)\widetilde{\bf q}\right \|_2\leq 10^{-14}\equiv tol_{\bf sys}
\end{equation}
as the stopping criterion for the Toeplitz systems. This algorithm mimics solving the Toeplitz systems ``exactly" via an iterative solver.
Let ${\bf y}(t)$ be the ``exact" solution and ${\bf y}_m(t)$ be the approximate solutions obtained from Algorithm 1 or Algorithm 2,
then we define
\begin{equation}\label{eqn402}
{\bf Error}=\frac{\|{\bf y}(t)-{\bf y}_m(t)\|_{2}}{\|{\bf y}(t)\|_{2}}
\end{equation}
as the relative error of the approximation ${\bf y}_m(t)$. Except for Example 1, the ``exact" solution ${\bf y}(t)$ is calculated by using the MATLAB bulit-in function {\it expm.m}.
In the tables below, we denote by {\bf CPU} the CPU time in seconds. We choose the vector ${\bf v}=[1,1,\ldots,1]^{\rm T}$ for all the numerical experiments in this section.

{\bf Example 1.}~~In this example, we aim to show the effectiveness of our inexact strategy (\ref{eqn5.21}), as well as the superiority of Algorithm 2 over Algorithm 1.
The Toeplitz matrix $A$ is generated by the even function $\theta^{2}$ defined on $[-\pi,\pi]$.
We want to compute ${\bf y}(t)={\rm exp}(-tA){\bf v}$ with $t=1,tol_{\bf exp}=10^{-6}$ and $n=1\times 10^5,2\times 10^5,\ldots,5\times 10^5$, respectively. Since the size $n$ of the Toeplitz matrix is very large, the MATLAB build-in function {\it expm.m} is infeasible for this problem. As a compromise, we run Algorithm 1 with the convergence tolerance $tol_{\bf exp}=10^{-14}$ for the ``exact" solution ${\bf y}(t)$. Table 1 lists the numerical results.

We see from Table 1 that Algorithm 2 converges much faster than Algorithm 1 in practical calculations, and the inexact strategy is both efficient and reliable.
Thanks to (\ref{eqn5.21}), it is only necessary to solve the Toeplitz system in the accuracy of $\mathcal{O}(10^{-9})$ instead of $10^{-14}$. Furthermore, the approximate solutions computed from the two methods have the same accuracy in terms of {\bf Error}.

{\small
\begin{table}[!h]
\begin{center}
\def\temptablewidth{0.8\textwidth}
{\rule{\temptablewidth}{1pt}}
\begin{tabular*}{\temptablewidth}{@{\extracolsep{\fill}}lccccc}

$n$ &{\bf Algorithm}&${\bf tol}_{\bf sys}$& {\bf Error}& ${\bf CPU}$  \\\hline


$1\times 10^{5}$         &Algorithm 1&$1.000\times 10^{-14}$& $4.615\times10^{-7}  $ &$2.531$\\
                 &Algorithm 2&$1.239\times 10^{-9}$&  $4.615\times10^{-7}  $ &$1.297$\\\hline

$2\times 10^{5}$        &Algorithm 1&$1.000\times 10^{-14}$& $3.263\times10^{-7}  $ &$7.360$\\
                  &Algorithm 2&$1.239\times 10^{-9}$& $3.263\times10^{-7}  $ &$3.422$\\\hline

$3\times 10^{5}$        &Algorithm 1 &$1.000\times 10^{-14}$&$2.664\times10^{-7}  $ &$13.375$\\
                  &Algorithm 2&$1.239\times 10^{-9}$ &$2.664\times10^{-7}  $ &$6.031$\\\hline

$4\times 10^{5}$        &Algorithm 1&$1.000\times 10^{-14}$&$2.307\times10^{-7}  $ &$20.313$\\
                  &Algorithm 2&$1.239\times 10^{-9}$&$2.307\times10^{-7}  $ &$9.125$\\\hline

$5\times 10^{5}$        &Algorithm 1&$1.000\times 10^{-14}$ &$2.064\times10^{-7}  $ &$27.719$\\
                  &Algorithm 2&$1.239\times 10^{-9}$&$2.064\times10^{-7}  $ &$10.187$\\
 \end{tabular*}
 {\rule{\temptablewidth}{1pt}}\\
 \end{center}
 \begin{center}
  {\small {\rm ~{\bf Table 1, Example 1:}~A comparison of Algorithm 1 and Algorithm 2, $t=1$, $\gamma=1/10$ and $tol_{\bf exp}=10^{-6}$.}}
 \end{center}
 \end{table}
}

{\bf Example 2.}~~The aim of this example is two-fold. First, we show the effectiveness of Theorem \ref{Thm3.3}. Second, we illustrate that
our proposed 1-norm ``GSF condition number" (\ref{eqn3155}) is a good estimation to the 1-norm ``classical condition number" of a Toeplitz matrix.
For the first aim, we run Algorithm 2 with the stopping criterion $tol_{\bf exp}$ chosen as $10^{-2},10^{-4},\ldots,10^{-10}$, and try to show that
$$
\|{\bf r}^{real}-{\bf r}^{comp}\|_{2}=\mathcal{O}(tol_{\bf exp}).
$$
In order to compute the ``real" residual, we first
form the approximation ${\bf y}_m(t)=V_m{\bf u}_{m}(t)$ explicitly, and then compute ${\bf r}^{real}$ by (\ref{eqn320}). The convergence threshold $tol_{\bf sys}$ for the Toeplitz systems is determined by using (\ref{eqn5.21}).

There are two test problems in this example, both of which are from
\cite{LP}. The first test matrix is the non-Hermitian Toeplitz
matrix generated by the function $ f(\theta)=\theta^{2}+{\bf
i}\cdot\theta^{3}, {\bf i}=\sqrt{-1},~\theta\in[-\pi,\pi]. $ Notice
that $\rm{Re}(f)=\theta^{2}\geq 0$ is an even function, and
$\rm{Im}(f)=\theta^{3}$ is an odd function. Table 2 lists the
numerical results of Algorithm 2 for ${\rm exp}(-tA){\bf v}$ with
$t=1,~\gamma=1/10$, and $n=3000$. For the first test problem, we
have $\kappa_1^{GSF}(I+\gamma A)\approx 1.275\times 10^{2}$, which
is of medium sized. In the second test problem, we consider pricing
options for a single underlying asset in Merton's jump-diffusion
model {\rm\cite{LP,RM}}. As the real part of the eigenvalues of the
Toeplitz matrix are less equal to zero, we are interested in
computing ${\rm exp}(tA){\bf v}$ with $t>0$, for more details, see
Example 3 of \cite{LP}. Table 3 gives the numerical results of
Algorithm 2 with $t=1,~\gamma=1$, and $n=3000$. For this test
problem, we have $\kappa_1^{GSF}(I+\gamma A)\approx 6.296\times
10^{7}$, which is relatively large.

Two remarks are in order. First, we see that $\|{\bf r}^{real}-{\bf r}^{comp}\|_{2}$ and $tol_{\bf exp}$ are about in the
same order in all the cases. This illustrates the effectiveness of Theorem \ref{Thm3.3}, as well as the efficiency of the inexact strategy (\ref{eqn5.21}).
Second, we observe from Table 2 and Table 3 that, if the GSF condition
number $\kappa_1^{GSF}(I+\gamma A)$ is medium sized, one can solve
the Toeplitz systems (\ref{eqn203}) with a relatively low accuracy.
Otherwise, we have to solve them with a relatively high accuracy.
For instance, if we choose
$tol_{\bf exp}=10^{-6}$, one has to solve the
Toeplitz systems in an accuracy of $\mathcal{O}(10^{-9})$ for the
first test problem, while an accuracy of $\mathcal{O}(10^{-13})$ is
required for the second test problem.


{\small
\begin{table}[!h]
\begin{center}
\def\temptablewidth{0.8\textwidth}
{\rule{\temptablewidth}{1pt}}
\begin{tabular*}{\temptablewidth}{@{\extracolsep{\fill}}lccc}
 ${\bf tol}_{\bf exp}$&${\bf tol}_{\bf sys}$&${\|{\bf r}^{real}-{\bf r}^{comp}\|_{2}}$ &{\bf Error} \\\hline

$10^{-2}$&$1.010 \times 10^{-5}$  &$1.519 \times 10^{-3}$  &$2.679\times 10^{-4}$          \\

$10^{-4}$&$1.010 \times 10^{-7}$  &$5.352 \times 10^{-5}$    &$6.480 \times 10^{-6}$   \\

$10^{-6}$&$1.010 \times 10^{-9}$ &$4.839 \times 10^{-5}$   &$3.985 \times 10^{-6}$   \\

$10^{-8}$&$1.010 \times 10^{-11}$ &$1.679 \times 10^{-8}$  &$1.701\times 10^{-9}$    \\

$10^{-10}$&$1.010 \times 10^{-13}$ &$2.667 \times 10^{-9}$    &$2.607 \times 10^{-10}$    \\



 \end{tabular*}
 {\rule{\temptablewidth}{1pt}}\\
 \end{center}
 \begin{center}
  {\small {\rm ~{\bf Table 2, the 1st test problem of Example 2:}~Numerical results of Algorithm 2 with different $tol_{\bf exp}$ for computing $e^{-tA}{\bf v}$,  $t=1,~\gamma=1/10$,~$n=3000$;
   $\kappa_1^{GSF}(I+\gamma A)\approx 1.275\times 10^{2}$.}}
 \end{center}
 \end{table}
}


{\small
\begin{table}[!h]
\begin{center}
\def\temptablewidth{0.8\textwidth}
{\rule{\temptablewidth}{1pt}}
\begin{tabular*}{\temptablewidth}{@{\extracolsep{\fill}}lccc}
 ${\bf tol}_{\bf exp}$&${\bf tol}_{\bf sys}$&${\|{\bf r}^{real}-{\bf r}^{comp}\|_{2}}$ &{\bf Error} \\\hline

$10^{-2}$&$4.236 \times 10^{-9}$  &$6.958\times 10^{-2}$&$9.375\times 10^{-5}$  \\

$10^{-4}$&$4.236  \times 10^{-11}$&$7.347 \times 10^{-4}$ &$6.817 \times 10^{-7}$  \\

$10^{-6}$&$4.236  \times 10^{-13}$&$3.320 \times 10^{-5}$&$3.056 \times 10^{-9}$   \\

$10^{-8}$&$4.236 \times 10^{-15}$ &$1.417 \times 10^{-7}$&$2.364\times 10^{-11}$   \\

$10^{-10}$&$4.236  \times 10^{-17}$ &$5.105 \times 10^{-11}$&$2.021 \times 10^{-11}$   \\



 \end{tabular*}
 {\rule{\temptablewidth}{1pt}}\\
 \end{center}
 \begin{center}
  {\small {\rm ~{\bf Table 3, the 2nd test problem of Example 2:}~Numerical results of Algorithm 2 with different $tol_{\bf exp}$ for computing $e^{tA}{\bf v}$, $t=1,~\gamma=1$,~$n=3000$;
  $\kappa_1^{GSF}(I+\gamma A)\approx 6.296\times 10^{7}$.}}
 \end{center}
 \end{table}
}
%
%
%
{\small
\begin{table}[!h]
\begin{center}
\def\temptablewidth{0.95\textwidth}
{\rule{\temptablewidth}{1pt}}
\begin{tabular*}{\temptablewidth}{@{\extracolsep{\fill}}lcccc}

 Test problem &$n$  & $\kappa_1(I+\gamma A)$& $\kappa_{1}^{\rm est}(I+\gamma A)$ &$\kappa_{1}^{GSF}(I+\gamma A)$\\\hline


                        &1000&$65.284$$(0.219)$  &$65.284(0.109)  $ &$79.037(0.031)$\\
1st test problem        &2000&$89.546(1.469)$  &$89.546(0.562)  $ &$1.071\times10^{2}(0.141)$\\
                       &3000&$1.073\times10^{2}(4.641)$  &$1.073\times10^{2}(1.547)  $ &$1.275\times10^{2}(0.265)$\\
                       &4000&$1.218\times 10^{2}(10.594)$  &$1.218\times 10^{2}(3.328)  $ &$1.442\times 10^{2}(0.500)$\\\hline

                       &1000&$2.436\times 10^{6}(0.219)$  &$2.436\times 10^{6}(0.094)  $ &$6.989\times 10^{6}(0.0620)$\\
2nd test problem       &2000&$9.736\times10^{6}(1.469)$  &$9.736\times10^{6}(0.531) $ &$2.797\times10^{7}(0.218)$\\
                       &3000&$2.190\times10^{7}(4.672)$  &$2.190\times10^{7}(1.531)  $ &$6.296\times10^{7}(0.453)$\\
                       &4000&$3.893\times 10^{7}(10.656)$  &$3.893\times 10^{7}(3.344)  $ &$1.119\times 10^{8}(0.828)$\\\hline

\end{tabular*}
 {\rule{\temptablewidth}{1pt}}\\
 \end{center}
 \begin{center}
  {\small {\rm {\bf Table 4, Example 2 }: The values of $\kappa_1^{GSF}(I+\gamma A)$, $\kappa_1^{\rm est}(I+\gamma A)$, $\kappa_1(I+\gamma A)$ and the CPU time in seconds for computing them (in brackets), $n=1000,2000,3000,4000$.}}
 \end{center}
 \end{table}
}
When $n=1000,2000,3000$ and 4000, we list in Table 4 the 1-norm GSF
condition number $\kappa_1^{GSF}(I+\gamma A)$( where we use
$\max{\|\bf fcol\|_{1}, \|\bf frow\|_{1}}$ instead of $\|I+\gamma
A\|_{1}$),
 the 1-norm classical condition number $\kappa_1(I+\gamma A)$ \big(evaluated by using the MATLAB command $cond(I+\gamma A,1)$\big) and its estimation $\kappa_1^{\rm est}(I+\gamma A)$ \big(evaluated by using the MATLAB command $condest(I+\gamma A)$\big); as well as the CPU time in seconds for solving them (in brackets). It is seen that the GSF condition number is about one to two times larger than the classical condition number, and the former is a good estimation to the latter.
Furthermore, the CPU time for $\kappa_1^{GSF}(I+\gamma A)$ is much less than that for $\kappa_1(I+\gamma A)$ and $\kappa_1^{\rm est}(I+\gamma A)$, especially when $n$ is large. Thus, the 1-norm GSF condition number is a competitive alternative to the classical condition number for Toeplitz matrices.

{\bf Example 3.}~In this example, we try to show that our new bounds (\ref{eqn3.4}) and (\ref{eqn325}) are sharper than (\ref{eqn.3.13}) and (\ref{eqn324}).
The test matrix is the ``gallery" matrix generated by the MATLAB command $A = gallery('parter',n)$ \cite{Higham}.
It is a Toeplitz matrix whose singular values are close to $\pi$. Let {\bf x} and {\bf y} be the ``exact" solutions of
the systems $(I +\gamma A){\bf x} = {\bf e}_{1}$ and $(I +\gamma A){\bf y} = {\bf e}_{n}$, respectively,
which are computed from running the preconditioned  (unrestarted) GMRES algorithm
with $tol_{\bf sys}=10^{-14}$. Then
we form  $\widetilde{\bf x}$ in the following way:
$$
{\bf f}=randn(n,1);\quad
{\bf f}={\bf f}/\|{\bf f}\|;\quad
\widetilde{\bf x}={\bf x}+\varepsilon\|{\bf x}\|\cdot {\bf f};
$$
where $randn(n,1)$ is a vector of length $n$ with normally
distributed random entries, and $\|\cdot\|$ is 1-norm for
(\ref{eqn3.4}) and (\ref{eqn325}), and 2-norm for (\ref{eqn.3.13})
and (\ref{eqn324}). The vector $\widetilde{\bf y}$ is formed in a
similar way. In this example, we choose
$\varepsilon=10^{-6},10^{-9},10^{-12}$ and $n=1000,2000,3000,4000$,
respectively. In order to show the sharpness of our results, we also
present the ``exact" absolute and relative errors
$\big\|T^{-1}-\widetilde{T}^{-1}\big\|_{2}$ and
$\frac{\left\|T^{-1}-\widetilde{T}^{-1}\right\|_{2}}{\left\|T^{-1}\right\|_2}$.
Tables 4 and 5 report the numerical results. It is seen that our
upper bounds are  sharper than those due to Gutknecht and Hochbruck,
especially when $n$ is large.

{\small
\begin{table}[!h]
\begin{center}
\def\temptablewidth{0.8\textwidth}
{\rule{\temptablewidth}{1pt}}
\begin{tabular*}{\temptablewidth}{@{\extracolsep{\fill}}lcccc}

 $\varepsilon$ &$n$&{ (\ref{eqn3.4}) }  & { (\ref{eqn.3.13})}& $\big\|T^{-1}-\widetilde{T}^{-1}\big\|_{2}$  \\\hline


                  &1000&$1.114\times10^{-5}$  &$3.358\times10^{-3}  $ &$3.238\times 10^{-6}$\\
$10^{-6}$         &2000&$1.217\times10^{-5}$  &$6.717\times10^{-3}  $ &$3.672\times 10^{-6}$\\
                 &3000&$1.281\times10^{-5}$  &$1.007\times10^{-2}  $ &$3.948\times 10^{-6}$\\
                 &4000&$1.329\times10^{-5}$  &$1.343\times10^{-2}  $ &$3.767\times 10^{-6}$\\\hline

                 &1000&$1.113\times10^{-8}$  &$3.358\times10^{-6}  $ &$3.182\times 10^{-9}$\\
$10^{-9}$        &2000&$1.218\times10^{-8}$  &$6.717\times10^{-6}  $ &$3.414\times 10^{-9}$\\
                  &3000&$1.282\times10^{-8}$  &$1.007\times10^{-5}  $ &$3.947\times 10^{-9}$\\
                   &4000&$1.328\times10^{-8}$  &$1.343\times10^{-5}  $ &$3.961\times 10^{-9}$\\\hline

                  &1000&$1.113\times10^{-11}$  &$3.358\times10^{-9}  $ &$3.385\times 10^{-12}$\\
$10^{-12}$        &2000&$1.217\times10^{-11}$  &$6.717\times10^{-9}  $ &$4.102\times 10^{-12}$\\
                  &3000&$1.282\times10^{-11}$  &$1.007\times 10^{-8}  $ &$3.930\times 10^{-12}$\\
                  &4000&$1.328\times10^{-11}$  &$1.343\times 10^{-8}  $ &$4.071\times 10^{-12}$\\

 \end{tabular*}
 {\rule{\temptablewidth}{1pt}}\\
 \end{center}
 \begin{center}
  {\small {\rm ~{\bf Table 5, Example 3:}~A comparison of the absolute error bounds (\ref{eqn3.4}) and (\ref{eqn.3.13}), $t=1,\gamma=1/10$, $\varepsilon=10^{-6},10^{-9},10^{-12}$ and $n=1000,2000,3000,4000$.}}
 \end{center}
 \end{table}
}

{\small
\begin{table}[!h]
\begin{center}
\def\temptablewidth{0.8\textwidth}
{\rule{\temptablewidth}{1pt}}
\begin{tabular*}{\temptablewidth}{@{\extracolsep{\fill}}lcccc}

 $\varepsilon$ &$n$&{ (\ref{eqn325}) }  & { (\ref{eqn324})}& $\frac{\left\|T^{-1}-\widetilde{T}^{-1}\right\|_{2}}{\left\|T^{-1}\right\|_2}$  \\\hline


                  &1000&$3.524\times10^{-4}$  &$4.813\times10^{-3}  $ &$2.932\times 10^{-6}$\\
$10^{-6}$         &2000&$5.447\times10^{-4}$  &$9.642\times10^{-3}  $ &$3.739\times 10^{-6}$\\
                 &3000&$7.018\times10^{-4}$  &$1.447\times10^{-2}  $ &$3.829\times 10^{-6}$\\
                 &4000&$8.403\times10^{-4}$  &$1.931\times10^{-2}  $ &$3.652\times 10^{-6}$\\\hline

                 &1000&$3.519\times10^{-7}$  &$4.813\times10^{-6}  $ &$3.626\times 10^{-9}$\\
$10^{-9}$        &2000&$5.443\times10^{-7}$  &$9.642\times10^{-6}  $ &$3.328\times 10^{-9}$\\
                  &3000&$7.018\times10^{-7}$  &$1.447\times10^{-5}  $ &$3.618\times 10^{-9}$\\
                   &4000&$8.403\times10^{-7}$  &$1.931\times10^{-5}  $ &$3.718\times 10^{-9}$\\\hline

                  &1000&$3.521\times10^{-10}$  &$4.813\times10^{-9}  $ &$3.427\times 10^{-12}$\\
$10^{-12}$        &2000&$5.443\times10^{-10}$  &$9.642\times10^{-9}  $ &$3.964\times 10^{-12}$\\
                  &3000&$7.017\times10^{-10}$  &$1.447\times 10^{-8}  $ &$4.893\times 10^{-12}$\\
                  &4000&$8.402\times10^{-10}$  &$1.931\times 10^{-8}  $ &$3.991\times 10^{-12}$\\

 \end{tabular*}
 {\rule{\temptablewidth}{1pt}}\\
 \end{center}
 \begin{center}
  {\small {\rm ~{\bf Table 6, Example 3:}~A comparison of the relative error bounds (\ref{eqn325}) and (\ref{eqn324}), $t=1,\gamma=1/10$, $\varepsilon=10^{-6},10^{-9},10^{-12}$ and $n=1000,2000,3000,4000$.}}
 \end{center}
 \end{table}
}
\section{Conclusion}
In this paper, we analyze and further develop an inexact shift-and-invert Arnoldi method for
the problem of numerical approximation to the product of Toeplitz matrix exponential with a vector. First, we
give an improved stability analysis on the Gohberg-Semencul formula (GSF) for the inverse of a Toeplitz matrix, and our result is independent of the size of the matrix in question. Moreover, we define the ``GSF condition number" of a Toeplitz matrix. An advantage is that we can evaluate the ``classical" condition number and the effective condition numbers of a Toeplitz matrix via solving Toeplitz systems, with no need to form the Toeplitz inverse explicitly. Second, we establish a relation between the error in approximating Toeplitz systems and the residual of its matrix exponential. Third, we provide a practical stopping criterion for the accuracy in approximating the Toeplitz systems in the inexact shift-and-invert Arnoldi algorithm for Toeplitz matrix exponential. It is shown that if the 1-norm ``GSF condition number" $\kappa_1^{GSF}(I+\gamma A)$ is medium sized, then the Toeplitz systems can be solved in a relatively low accuracy. 

\section*{Acknowledgments}
We would like to express our sincere thanks to Prof. Panayot Vassilevski and  two reviewers
for their invaluable comments and
constructive suggestions which greatly improve the presentation
of this paper.



\begin{thebibliography}{99}
{\small

\bibitem{MA} M. Abdou and A. Badr. {\rm On a method for solving an integral equation in the displacement contact problem}. Applied Mathematics and Computation 2002; {\bf 127}: 65--78.

\bibitem{RRA} M. Botchev, V. Grimm, and M. Hochbbruck. {\rm Residual, restarting and Richarson iteration for the matrix exponential}. SIAM Journal on Scientific Computing  2013;
{\bf 35}: A1376--A1397.

\bibitem{CJ} R. Chan and  X. Jin. {\it An Introduction to Iterative Toeplitz Solvers}. Society for Industrial and Applied Mathematics (SIAM): Philadelphia, PA, 2007.


\bibitem{CN} R. Chan and  M. Ng. {\it Conjugate gradient methods for Toeplitz systems}. SIAM Review, 1996; {\bf 38}: 427--482.

\bibitem{TChan} T. Chan. {\rm An optimal circulant preconditioner for Toeplitz systems}. SIAM Journal on Scientific and Statistical Computing 1988; {\bf 9}: 766--771.

\bibitem{CMSW} A. Cline, C. Moler, G.W. Stewart, and J.H. Wilkinson. {\rm An estimate for the condition number of a matrix}.
SIAM Journal on Numerical Analysis 1979; {\bf  16}: 368--375.

\bibitem{Our} T. Feng, G. Wu, and T. Xu. {\rm An inexact shift-and-invert Arnoldi algorithm for large non-Hermitian generalized Toeplitz eigenproblems}. submitted, 2014.

\bibitem{Frommer} A. Frommer, K. Kahl, T. Lippert, and H. Rittich.
{\rm 2-norm error bounds and estimates for Lanczos approximations to linear systems and rational matrix functions}.
SIAM Journal on Matrix Analysis and Applications 2013; {\bf 34}: 1046--1065.

\bibitem{GSF} I. Gohberg and  A. Semencul. {\rm On the inversion of finite Toeplitz matrices and their continuous analogs}. Mat. Issled.,1972; {\bf  2}: 201--233.

\bibitem{GHK} I. Gohberg, M. Hanke, and I. Koltracht. {\rm Fast preconditioned conjugate gradient algorithms for Wiener-Hopf
integral equations}. SIAM Journal on Numerical Analysis 1994; {\bf 31}: 429--443.

\bibitem{GV} G.H. Golub and  C.F. Van Loan. {\it Matrix Computations}.
4th edition, John Hopkins University Press: Baltimore, MD, 2013.

\bibitem{MHG} M. Gutknecht and  M. Hochbruck. {\rm The stability of inversion formulas for Toeplitz matrices}. Linear Algebra and its Applications 1995;  {\bf 223/224}: 307--324.

\bibitem{Higham1} N.J. Higham. {\it Functions of Matrices: Theory and Computation}. Society for Industrial and Applied Mathematics (SIAM): Philadelphia, PA, 2008.

\bibitem{Higham} N.J. Higham. {\it Accuracy and Stability of Numerical Algorithms}. 2nd edition, Society for Industrial and Applied Mathematics (SIAM): Philadelphia, PA, 2002.

\bibitem{LP} S. Lee, H. Pang, and H. Sun. {\rm Shift-invert Arnoldi approximation to the Toeplitz matrix exponential}. SIAM Journal on Scientific Computing 2010; {\bf 32}: 774--792.

\bibitem{Li} Z. Li, H. Huang, Y. Wei and A. Chen. {\it Effective Condition Number for Numerical Partial Differential Equations}.  Science
Press: Beijing, 2013 and Alpha Science International Ltd.: Oxford, 2014.




\bibitem{RM} R. Merton. {\rm Option pricing when underlying stock returns are discontinuous}. Journal of Financial Economics 1976; {\bf  3}: 125--144.

\bibitem{ML} C. Moler and  C.F. Van Loan. {\rm Nineteen dubious ways to compute the exponential of a matrix,
twenty-five years later}. SIAM Review 2003, {\bf  45}: 3--49.

\bibitem{MN} I. Moret and  P. Novati. {\rm RD-rational approximations of the matrix exponential}. BIT 2004;  {\bf 44}: 595--615.

\bibitem{Ng2} M. Ng. {\it Iterative Methods for Toeplitz Systems}. Oxford University Press: New York,  2004.

\bibitem{Pang} H. Pang and  H. Sun. {\rm Shift-Invert Lanczos method for the symmetric positive semidefinite Toeplitz matrix exponential}. Numerical Linear Algebra with Applications 2011; {\bf 18}: 603--614.

\bibitem{Rice} J. Rice. {\it Matrix Computations and Mathematical Software}. McGraw-Hill Book Company:  New York, 1981.

\bibitem{Saad} Y. Saad. {\rm Analysis of some Krylov subspace approximations to the matrix exponential operator}.
 SIAM Journal on Numerical Analysis 1992; {\bf  29}: 209--228.

\bibitem{Saad2} Y. Saad.
{\it Iterative Methods for Sparse Linear Systems}. 2nd edition, Society for Industrial and Applied Mathematics
(SIAM): Philadelphia, PA, 2003.

\bibitem{DAM} D. Tangman, A. Gopaul, and M. Bhuruth. {\rm Exponential time integration and Chebychev discretisation schemes for fast pricing of options}. Applied Numerical Mathematics 2008;  {\bf 58}: 1309--1319.

%

\bibitem{JH} J. van den Eshof and  M. Hochbruck. {\rm Preconditioning Lanczos approximations to the matrix exponential}. SIAM Journal on Scientific Computing 2006; {\bf  27}: 1438--1457.

}
\end{thebibliography}
\end{document}